\documentclass[11pt]{article}
\usepackage{amsmath}
\usepackage{graphicx}
\usepackage{enumerate}
\usepackage[numbers]{natbib}
\usepackage{url} 
\usepackage{epstopdf}
\usepackage{epsfig}
\usepackage{multirow}
\usepackage{graphics}
\usepackage{subfigure}
\usepackage{caption}
\usepackage{float}
\usepackage{color}
\usepackage{array}
\usepackage{booktabs}
\usepackage{subfig}
\usepackage{lscape}
\usepackage{booktabs}

\newtheorem{defn}{Definition}
\newtheorem{lemma}{Lemma}
\newtheorem{prop}{Proposition}

\begin{document}

	\def\spacingset#1{\renewcommand{\baselinestretch}%
		{#1}\small\normalsize} \spacingset{1}

	
	\title{\bf Simpson's Paradox with Any Given Number of Factors}
	
		\author{Guisheng Dai$^{1}$ and Weizhen Wang$^{2}$
		\thanks{{Corresponding author. Department of Mathematics and Statistics, Wright State University, Dayton, Ohio 45435, USA} \mbox{\emph{E-mail~address}: weizhen.wang@wright.edu}}
		\\
		\footnotesize
		$^{1}$\textit{School of Mathematics, Beijing Normal University, Beijing 100875, P. R. China} \vspace{0.0cm}
		\\\footnotesize
		$^{2}$\textit{Department of Mathematics and Statistics, Wright State University, Dayton, OH 45435, U.S.A} \vspace{-0.5cm}
		\\\footnotesize
	}
	\maketitle
\bigskip
\begin{abstract} 
	
	Simpson’s Paradox is a well-known phenomenon in statistical science, where the relationship
	between the response variable $X$ and a certain explanatory factor of interest $A$ reverses when an additional factor $B_1$ is considered. This paper explores
	the extension of Simpson’s Paradox to any given number $n$ of factors, referred to as the $n$-factor
	Simpson’s Paradox. We first provide a rigorous definition of the $n$-factor Simpson’s Paradox, then
	demonstrate the existence of a probability distribution through a geometric construction. Specifically, we show that for any positive
	integer $n$, it is possible to construct a probability distribution in which the conclusion about the effect of 
	$A$ on $X$ reverses each time an additional factor $B_i$ is introduced for $i=1,...,n$. A detailed example for $n = 3$ illustrates the construction. Our results highlight that, contrary to the intuition that more data
	leads to more accurate inferences, the inclusion of additional factors can repeatedly reverse
	conclusions, emphasizing the complexity of statistical inference in the presence of multiple
	confounding variables. 
	
\end{abstract}

\noindent%
{\it Keywords: Angle of vector; Explanatory variable; Sum of vectors; Tangent function.}

\vfill
\newpage
\spacingset{1.5} 

\section{Introduction}
\setcounter{equation}{0} 

Simpson's Paradox is an interesting phenomenon in which the conclusion of the effect of a certain explanatory variable  of interest $A$ on the response variable $X$ may reverse when another explanatory variable $B_1$ is observed. We validate this well-known phenomenon in elementary plane geometry. Also, when up to $n$ explanatory variables $B_1,...,B_n$ step in the study sequentially one by one, is it possible that the conclusion reverse each time when a new explanatory variable is observed?  A positive answer would indicate that more data information may lead to a sequence of inconsistent conclusions. This may further raise a concern of the usefulness of big data when establishing a treatment effect. 

Simpson (1951) first described this phenomenon in a technical paper, but Pearson (1899) and Yule (1903) had mentioned similar effects earlier. The name Simpson's Paradox was introduced in Blyth (1972). Bickel, Hammel, and O’Connell (1975) addressed the sex bias in graduate admissions through Simpson's Paradox. 
For simplicity, we consider the problem when all explanatory variables in the study assume two levels, as an extension to more than two levels is similar. Our approach to constructing Simpson's Paradox is to find two groups of two lines (i.e., four lines in total) with appropriate angles using  two initial given lines. This is accessible to students with knowledge of conditional probability and plane geometry.

Let $X$ be a binary random variable assuming two values: $X_1$ (Success) or $X_0$ (Failure).
Let $A$ be the factor of interest that assumes two levels: $A_1$ (Treatment 1) and $A_0$ (Treatment 0). Let $n$ be a given positive integer and let $B_i$ be an additional factor of two levels: $B_{i,1}$ and $B_{i,0}$ for $i=1,...,n$. In the study we observe $(X,A)$, as well as $(B_1,...,B_m)$ for some  integer $m\in [0,n]$.  

For example, $X$ is the status of a patient after seeing a doctor in a hospital: $X_1$ (cured) and $X_0$ (not cured); 
$A$ is the factor of two hospitals: $A_1$ (a local clinic) and $A_0$ (a national hospital); $B_1$ is the health condition of patient: $B_{1,1}$ (severe) and  $B_{1,0}$ (not severe);  $B_2$ is the location of patient: $B_{2,1}$ (city) and  $B_{2,0}$ (rural);  $B_3$ is the income level of patient: $B_{3,1}$ (high) and  $B_{3,0}$ (not high); etc. 
The goal is to determine the effect of hospital on the curing rate  utilizing information from up to $n$ possible factors $B_i$'s sequentially. 

For two events $C$ and $D$, write the intersection of $C$ and $D$ as $CD$ rather than $C\cap D$.

When $A$ is available only without any $B_i$'s, we
compare one pair of probabilities
$$P(X_1|A_1) \,\ vs. \,\ P(X_1|A_0).$$
If we find $P(X_1|A_1) > P(X_1|A_0)$, it is natural to conclude that $A_1$ is better than $A_0$. 
When $A$ and $B_1$ are available without $B_i$'s for $i\geq 2$, we compare two pairs of probabilities
$$P(X_1|A_1B_{1,1}) \,\ vs. \,\ P(X_1|A_0B_{1,1})\,\ and \,\ P(X_1|A_1B_{1,0}) \,\ vs. \,\ P(X_1|A_0B_{1,0}).$$
If we find both $P(X_1|A_1B_{1,1})< P(X_1|A_0B_{1,1})$ and $P(X_1|A_1B_{1,0})< P(X_1|A_0B_{1,0})$, then we intend to conclude that $A_1$ is worse than $A_0$, which reverses the conclusion when only $A$ is considered. i.e.,  
Simpson's Paradox occurs. 

An example was given in  Chang et al. (1986), where they compared the success rates (success $=X_1$, ) of two treatments, $A_1$ and $A_0$, for kidney stones in 700 patients, with 350 patients assigned to each treatment.  Initially, $A_1$ had a higher success rate than $A_0$; however, when the stone size ($B_1$) was taken into account, this conclusion was reversed.  
\begin{center}
	\begin{tabular}{llcl}
		\hline
		& Treatment $A_1$ & Association & Treatment $A_0$ \\
		\hline
		No size          & $P(X_1|A_1)=83\% (289/350)$ &$>$& $P(X_1|A_0)=78\% (273/350)$ \\
		\hline
		Small stone ($B_{1,1}$) & $P(X_1|A_1B_{1,1})=87\% (234/270)$  &$<$& $P(X_1|A_1B_{1,0})=93\% (81/87)$ \\
		Large stone ($B_{1,0}$) & $P(X_1|A_1B_{0,1})=69\% (55/80)$ & $<$& $P(X_1|A_1B_{0,0})=73\% (192/263)$ \\
		\hline
	\end{tabular}
\end{center} 

\vspace{0.1in}

In general, for any positive integer $m\leq n$, when $A$, $B_1,...,B_m$ are available,
let $$S_m=\{s_m=(i_1,...,i_m): i_j=0, 1, j\in [1,m]\}$$ be the set of all possible outcomes of $(B_1,...,B_m)$.  Define a point in $S_m$
$$\bar{B}_{s_m}=\{B_1=i_1,...,B_m=i_m\}.$$ 
For example, when $m=1$, then $s_m=1$ or 0, and  $B_{1,1}=\bar{B}_{s_m}$ for $s_m=1$ and  $B_{1,0}=\bar{B}_{s_m}$ for $s_m=0$. 
For simplicity, when $m=0$, denote $\bar{B}_{s_0}=S$, the sample space.
We compare
$$P(X_1|A_1\bar{B}_{s_m}) \,\ vs. \,\ P(X_1|A_0\bar{B}_{s_m}),\,\ \forall s_m\in S_m, \forall \,\ 0\leq m\leq n.$$

We wish to construct a probability distribution so that
\begin{equation}\label{sim-0}
	P(X_1|A_1)=P(X_1|A_1\bar{B}_{s_0}) > P(X_1|A_0)=P(X_1|A_0\bar{B}_{s_0});
\end{equation}
\begin{equation}\label{sim-1}
	P(X_1|A_1\bar{B}_{s_1}) < P(X_1|A_0\bar{B}_{s_1}),\,\ \forall s_1\in S_1;
\end{equation}
\begin{equation}\label{sim-2}
P(X_1|A_1\bar{B}_{s_2}) > P(X_1|A_0\bar{B}_{s_2}),\,\ \forall s_2\in S_2;
\end{equation}
$$...$$
\begin{equation}\label{sim-m}
\left\{	\begin{split}
		&P(X_1|A_1\bar{B}_{s_m}) > P(X_1|A_0\bar{B}_{s_m}),\,\ \forall s_m\in S_m, \,\ \mbox{if $m$ is even};\\ 	
		&P(X_1|A_1\bar{B}_{s_m}) < P(X_1|A_0\bar{B}_{s_m}),\,\ \forall s_m\in S_m, \,\ \mbox{if $m$ is odd};
	\end{split}
\right.
\end{equation}
$$...$$
\begin{equation}\label{sim-n}
\left\{	\begin{split}
		&P(X_1|A_1\bar{B}_{s_n}) > P(X_1|A_0\bar{B}_{s_n}),\,\ \forall s_n\in S_n, \,\ \mbox{if $n$ is even};\\ 	
		&P(X_1|A_1\bar{B}_{s_n}) < P(X_1|A_0\bar{B}_{s_n}),\,\ \forall s_n\in S_n, \,\ \mbox{if $n$ is odd}.
	\end{split}
\right.
\end{equation}
i.e., the direction of inequality does not change for all $s_m$'s when $m$ is fixed but switches to the opposite as $m$ increases by one. This means that if none of the $B_i$'s are available, we conclude the treatment of $A_1$ better than the treatment of $A_0$ due to (\ref{sim-0}); if only $B_1$ is available, we conclude  $A_1$  worse than $A_0$ due to (\ref{sim-1}); if only $B_1$ and $B_2$ are available, we conclude $A_1$ better than $A_0$ again due to (\ref{sim-2}), etc. Hence, the effect of $A$ on $X$ cannot be reliably inferred as the conclusion is flipped over and over as the information in an additional $B_m$ is collected each time. 

When considering $
m=0$ and 
$m=1$ only, i.e., (\ref{sim-0}) and (\ref{sim-1}), the direction of the inequalities changes once, and this represents the classic Simpson’s Paradox. When considering 
$m=0,1,…,n$ for a general 
$n>1$, i.e., (\ref{sim-0}) through (\ref{sim-n}), the direction of the inequalities changes 
$n$ times. This phenomenon may occur in observational studies where the assignment of levels in each $B_i$
over experimental units (e.g., patients) is not random and may depend on 
$A$. 
We propose the following concept.
\begin{defn}
A probability distribution on the random vector $(X,A,B_1,...,B_n)$ that satisfies the relationships (\ref{sim-0}) though (\ref{sim-n}) is called the $n$-factor Simpson's Paradox. When $n=1$, it is simply called  the Simpson’s Paradox.
\end{defn}


Simpson’s paradox is not difficult to characterize mathematically, but is challenging to express
intuitively. Lindley and Novick (1981) proved the existence of Simpson’s Paradox with one explanatory
variable by using probability theory. Good and Mittal (1987) explained why Simpson’s Paradox
occurs using simple measures from 2 × 2 contingency tables and suggested ways to prevent it.
Kolick (2001) used planar geometric methods to intuitively explain the reasons behind Simpson’s
paradox with one explanatory variable.

For any four given positive real numbers $a_0,b_0,c_0$ and $d_0$ in $[0,1]$ satisfying \begin{equation}\label{sim-00}
	{a_0\over a_0+b_0}>(\mbox{or} <){c_0\over c_0+d_0},
\end{equation} we aim to offer a geometric construction of the $n$-factor Simpson's Paradox satisfying $$P(X_1|A_1)={a_0\over a_0+b_0}\,\ \mbox{and}\,\ P(X_1|A_0)={c_0\over c_0+d_0}.$$ 
\section{A geometric construction of the $n$-factor Simpson's Paradox}

The $n$-factor Simpson's Paradox is constructed by induction starting from the Simpson's Paradox and consists of $2^n-1$  Simpson's Paradoxes. For example,  a $2$-factor Simpson's Paradox may consist of the following $3$ different Simpson's Paradoxes: 
\begin{equation}\label{sim-2-all}
	\begin{split}
		&\mbox{i) $P(X_1|A_1) > P(X_1|A_0)$ vs. $P(X_1|A_1B_{1,1})<P(X_1|A_0B_{1,1})$ and $P(X_1|A_1B_{1,0})<P(X_1|A_0B_{1,0})$};\\ 
		& \mbox{ii)  $P(X_1|A_1B_{1,1})<P(X_1|A_0B_{1,1})$ vs.  $P(X_1|A_1B_{1,1}B_{2,i})>P(X_1|A_0B_{1,1}B_{2,i}$) for $i=1,0$};\\
		& \mbox{iii)  $P(X_1|A_1B_{1,0})<P(X_1|A_0B_{1,0})$ vs.  $P(X_1|A_1B_{1,0}B_{2,i})>P(X_1|A_0B_{1,0}B_{2,i}$) for $i=1,0$}.\\
		\end{split}
\end{equation}
 Hence,  it is enough to describe a geometric construction of the Simpson's Paradox. i.e., construct a distribution  from  any four given positive real numbers $a_0,b_0,c_0$ and $d_0$ satisfying (\ref{sim-00}) so that  
\begin{equation}\label{sim-0-1}
	P(X_1|A_1) > P(X_1|A_0);
\end{equation}
\begin{equation}\label{sim-1-1}
	P(X_1|A_1B_{1,1}) < P(X_1|A_0B_{1,1}), P(X_1|A_1B_{1,0}) < P(X_1|A_0B_{1,0}).
\end{equation}

First, decompose the vector of $(b_0,a_0)$ in the first quadrant of the $R^2$ plane as the sum of any two vectors $(b_1,a_1)$ and $(b_2,a_2)$ also in the first quadrant. i.e.,
\begin{equation}\label{decom-1}
	(b_0,a_0)=(b_1,a_1)+(b_2,a_2).
\end{equation}
 Hence, $a_i$'s and $b_i$'s are all positive. 
Let $\theta_i$ be the angle of vector $(b_i,a_i)$ for $i=0,1,2$. So, $\theta_i\in (0,\pi/2)$ and $\tan(\theta_i)=a_i/b_i$. 
In the setting of the Simpson's Paradox, without loss of generality, assume $a_0+b_0=1$. Indeed, we have
\begin{equation}\label{dist-1}
\left\{ 
\begin{split}
	& a_0=P(X_1|A_1), b_0=P(X_0|A_1), \\
	& a_1=P(X_1B_{1,1}|A_1), a_2=P(X_1B_{1,0}|A_1); b_1=P(X_0B_{1,1}|A_1),b_2=P(X_0B_{1,0}|A_1).
\end{split}
\right.
\end{equation} 
Therefore,
$${a_1\over a_1+b_1}=P(X_1|A_1B_{1,1}), {a_2\over a_2+b_2}=P(X_1|A_1B_{1,0}).$$

Second, repeat the process in the last paragraph to decompose vector $(d_0,c_0)$ as
\begin{equation}\label{decom-2}
	(d_0,c_0)=(d_1,c_1)+(d_2,c_2),
\end{equation}
 and 
let $\eta_i$ be the angle of vector $(d_i,c_i)$ for $i=0,1,2$.
In the setting of the Simpson's Paradox,  without loss of generality, also assume $c_0+d_0=1$. Indeed, we have
\begin{equation}\label{dist-2}
	\left\{ 
	\begin{split}
		&c_0=P(X_1|A_0), d_0=P(X_0|A_0), \\ &c_1=P(X_1B_{1,1}|A_0),c_2=P(X_1B_{1,0}|A_0), d_1=P(X_0B_{1,1}|A_0),  d_2=P(X_0B_{1,0}|A_0).
	\end{split}
\right.
\end{equation}
Therefore,
$${c_1\over c_1+d_1}=P(X_1|A_0B_{1,1}), {c_2\over c_2+d_2}=P(X_0|A_0B_{1,0}),$$
and (\ref{sim-1-1}) is equivalent to
$${a_1\over a_1+b_1}<{c_1\over c_1+d_1}, {a_2\over a_2+b_2}<{c_2\over c_2+d_2}. $$

\begin{lemma}\label{lem-equi} For any four positive constants $x_1,x_2,y_1,$ and $y_2$, let $\theta_x$ and $\theta_y$ be the angles of vectors $(x_2,x_1)$ and $(y_2,y_1)$ in the $R^2$ plane, respectively. Then
the following three are equivalent: 

i) 	${x_1\over x_1+x_2}>{y_1\over y_1+y_2}$; 

ii) $\tan(\theta_x)>\tan(\theta_y)$;

iii) $\theta_x>\theta_y$.
\end{lemma}

{\bf Proof}. Note $\tan(\theta_x)=x_1/x_2$, $\tan(\theta_y)=y_1/y_2$, and both $\theta_x$ and $\theta_y$ are in $ (0,\pi/2)$. The equivalence between i) and ii) follows
$$\tan(\theta_x)>\tan(\eta_y) \Leftrightarrow {x_1\over x_2}>{y_1\over y_2} \Leftrightarrow {x_1\over x_1+x_2}>{y_1\over y_1+y_2}.$$ The equivalence between ii) and iii) is due to the fact that function $\tan$ is strictly increasing over interval $(0,\pi/2)$. \fbox{}

If applying Lemma~\ref{lem-equi} to $(b_0,a_0)$ with angle $\theta_0$ and $(d_0,c_0)$ with angle $\eta_0$, it is clear that (\ref{sim-0-1}) is equivalent to 
$\theta_0>\eta_0$. If applying Lemma~\ref{lem-equi} to $(b_i,a_i)$ with angle $\theta_i$ and $(d_i,c_i)$ with angle $\eta_i$ for $i=1,2$, then (\ref{sim-1-1}) is equivalent to 
$\theta_i<\eta_i$ for $i=1,2$.
To achieve (\ref{sim-1-1}), we have the lemma below.

\begin{lemma}\label{lem-simpson-1} For any four positive constants $a_0,b_0,c_0,$ and $d_0$ satisfying (\ref{sim-00}), which  is equivalent to $\pi/2>\theta_0>\eta_0>0$, there exists a decomposition of (\ref{decom-1}) and (\ref{decom-2}) such that 
	\begin{equation}\label{sim-g-1}
			0<\theta_i<\eta_i<\pi/2, \,\  i=1,2.
	\end{equation}
	The probability distribution is given in (\ref{dist-1}) and (\ref{dist-2}).
\end{lemma}

We see that the direction of the inequalities reverses as factor $B_1$ steps in. i.e., the Simpson's Paradox occurs. A parallel result of another type of Simpson's Paradox is stated below.

\begin{lemma}\label{lem-simpson-1-1} For any four positive constants $a_0,b_0,c_0$, and $d_0$ satisfying  $0<\theta_0<\eta_0<\pi/2$, there exists a decomposition of (\ref{decom-1}) and (\ref{decom-2}) such that $$\pi/2>\theta_i>\eta_i>0, \,\  i=1,2.$$
	The probability distribution is given in (\ref{dist-1}) and (\ref{dist-2}).
\end{lemma}

\begin{figure}[htbp]
 	\centering
 	\vspace{-2mm}
 	\subfigure[Definition of $\theta_1$ and $\eta_1$]{
 		\begin{minipage}[b]{.3\linewidth}
 			\centering
 			\includegraphics[scale=0.32]{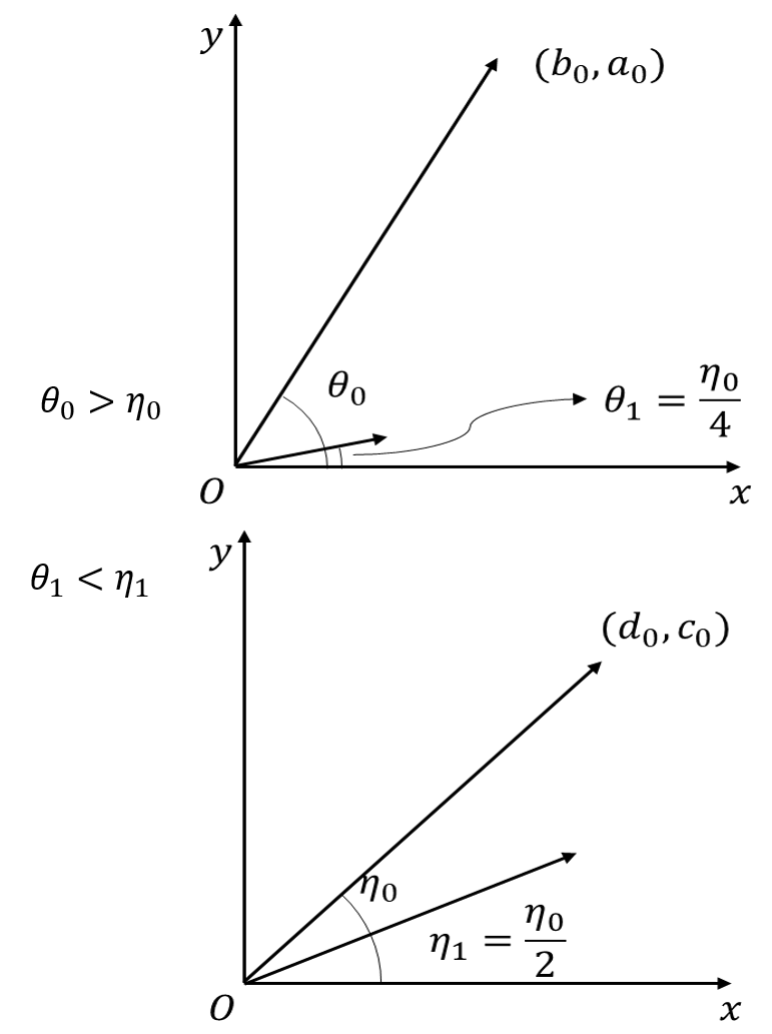}
 		\end{minipage}
 	}
 	\subfigure[Definition of $\theta_2$ and $\eta_2$]{
 		\begin{minipage}[b]{.3\linewidth}
 			\centering
 			\includegraphics[scale=0.32]{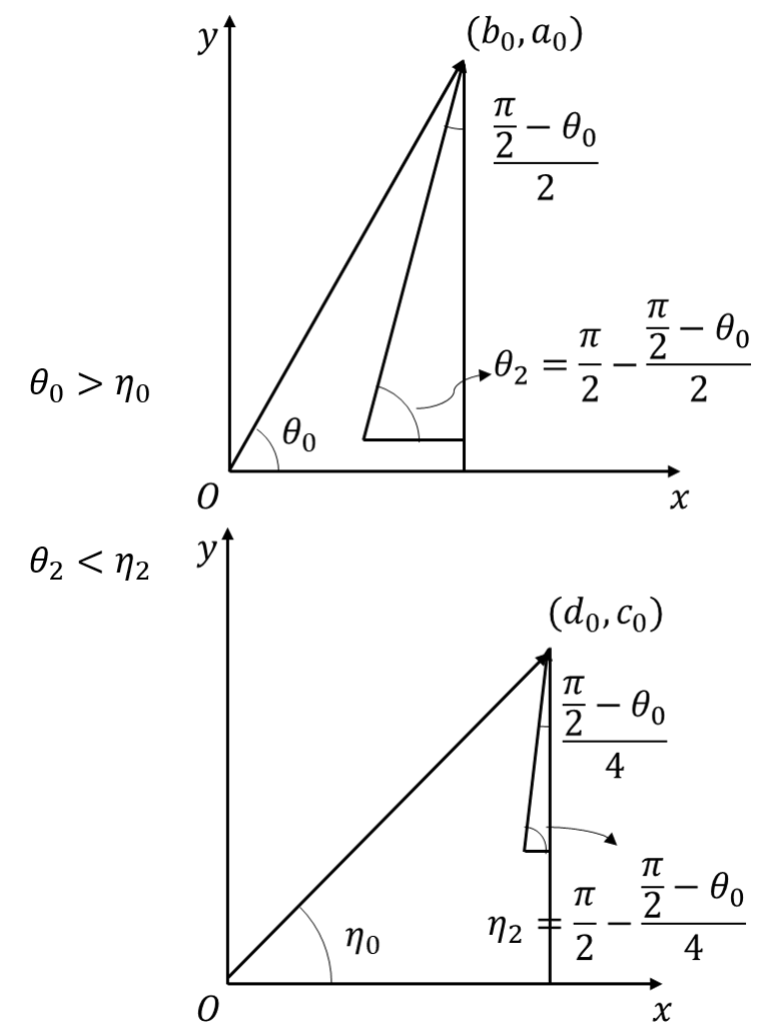}
 		\end{minipage}
 	}
 	\subfigure[Result of Lemma~\ref{lem-simpson-1}]{
 		\begin{minipage}[b]{.3\linewidth}
 			\centering
 			\includegraphics[scale=0.32]{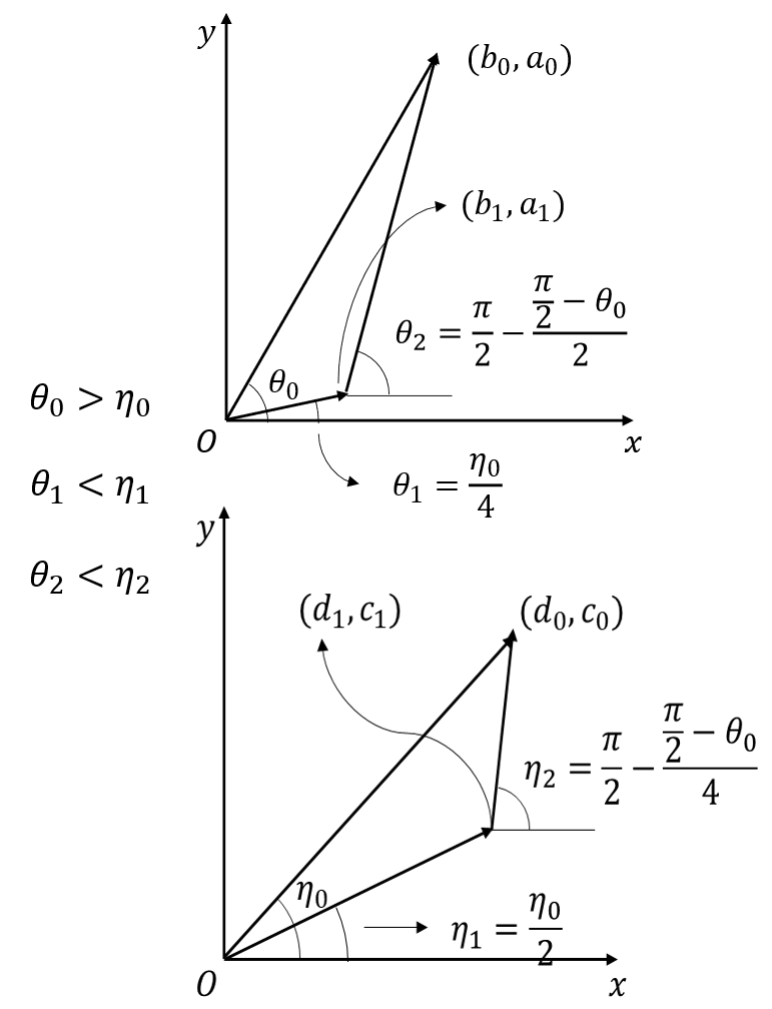}
 		\end{minipage}
 	}
 	\vspace{-2mm}
 	\subfigure[Definition of $\theta_1$ and $\eta_1$]{
 		\begin{minipage}[b]{.3\linewidth}
 			\centering
 			\includegraphics[scale=0.32]{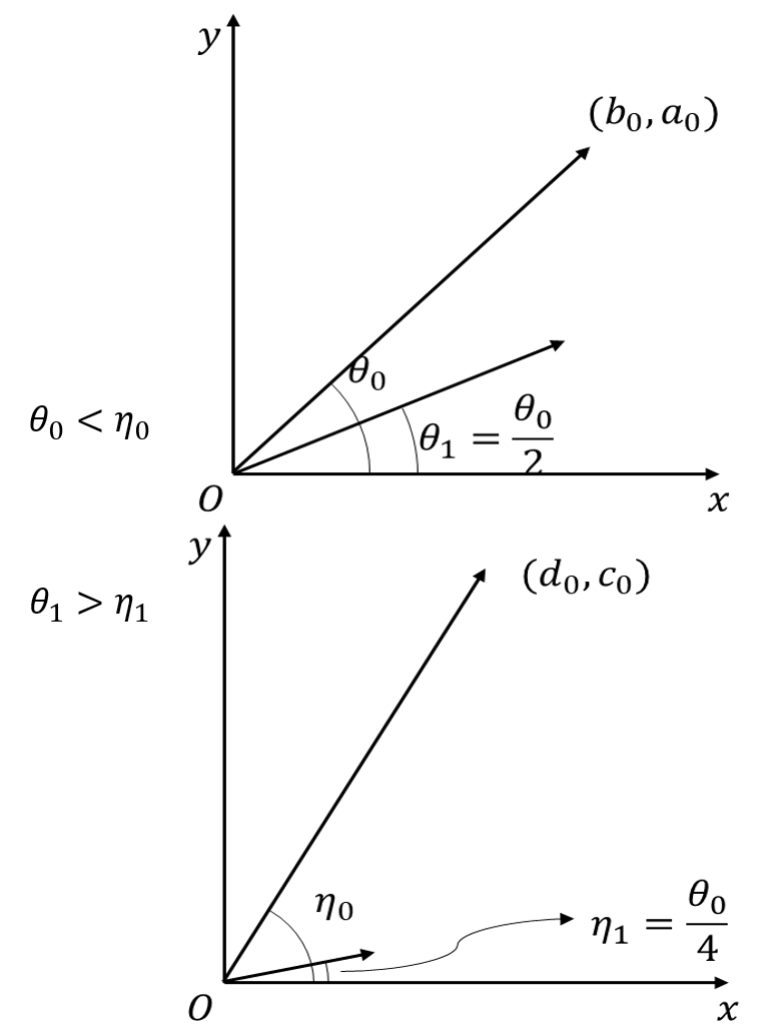}
 		\end{minipage}
 	}
 	\subfigure[Definition of $\theta_2$ and $\eta_2$]{
 		\begin{minipage}[b]{.3\linewidth}
 			\centering
 			\includegraphics[scale=0.32]{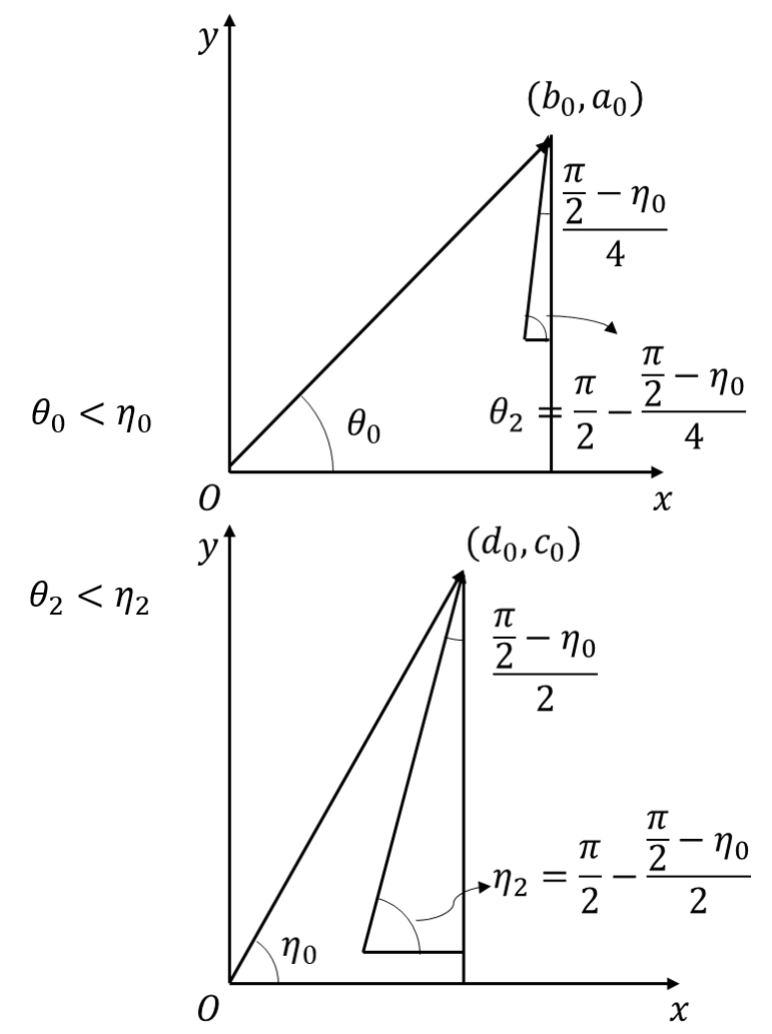}
 		\end{minipage}
 	}
 	\subfigure[Result of Lemma~\ref{lem-simpson-1-1}]{
 		\begin{minipage}[b]{.3\linewidth}
 			\centering
 			\includegraphics[scale=0.32]{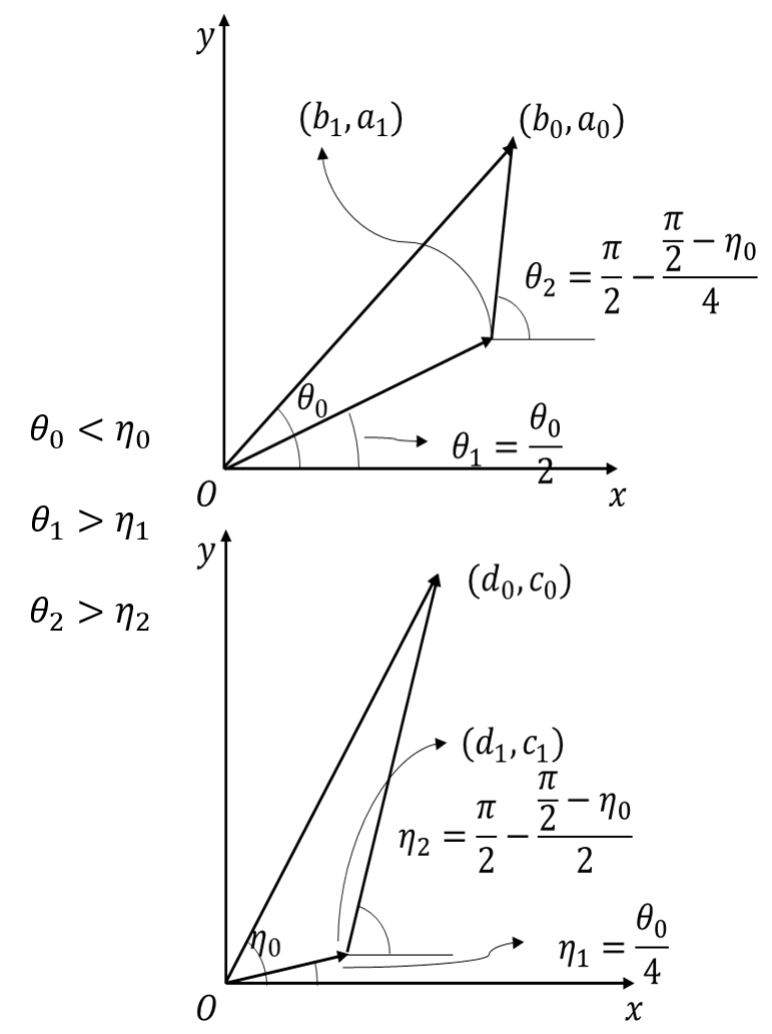}
 		\end{minipage}
 	}
 	\caption{The construction of the Simpson's Paradox. The top six graphs illustrate  Lemma~\ref{lem-simpson-1} and the bottom six are for Lemma~\ref{lem-simpson-1-1}.}
 	\label{figure-simpson}
 \end{figure}

{\bf Proof}. We first prove Lemma~\ref{lem-simpson-1}. 
Due to the decomposition (\ref{decom-1}) and the fact that all vectors of $(b_i,a_i)$ and $(d_i,c_i)$ are in the first quadrant in $R^2$, we have $\theta_1 \in(0, \theta_0)$, $\theta_2\in (\theta_0,\pi/2)$; $\eta_1 \in(0, \eta_0)$, $\eta_2\in (\eta_0,\pi/2)$. There exist many choices of $\theta_i$ and $\eta_i$ satisfying (\ref{sim-g-1}). For example, pick 
\begin{equation}\label{sim-geo-1}
	\left\{
	\begin{split}
&		\theta_1={\eta_0\over 4},\,\   \eta_1={\eta_0\over 2};\\
&		\theta_2={\pi\over 2}-{({\pi\over 2}-\theta_0)\over 2}={\pi\over 4}+{\theta_0\over 2}, \,\ \eta_2={\pi\over 2}-{{\pi\over 2}-\theta_0\over 4}={3\pi\over 8}+{\theta_0\over 4}.
	\end{split}
	\right.
\end{equation}
 So, (\ref{sim-g-1}) is true. See this construction in Figure~\ref{figure-simpson}--a,b,c. The proof of Lemma~\ref{lem-simpson-1} is complete.
 
 Lemma~\ref{lem-simpson-1-1} is proved similarly by picking
 \begin{equation}\label{sim-geo-2}
 	\left\{
 	\begin{split}
 		&	\theta_1={\theta_0\over 2}, \,\ \eta_1={\theta_0\over 4};\\
 		&	\theta_2
 		={3\pi\over 8}+{\eta_0\over 4},\,\ \eta_2
 		={\pi\over 4}+{\eta_0\over 2}.
 	\end{split}
 	\right.
 \end{equation}
 This construction is also displayed in Figure~\ref{figure-simpson}-- d,e,f.
 \fbox{}
 
To implement the decomposition in (\ref{sim-geo-1} and (\ref{sim-geo-2}) numerically, we provide the values of $a_i,b_i,c_i$ and $d_i$ which correspond to $\theta_i$ and $\eta_i$ for $i=1,2$ below for  given $a_0,b_0,c_0$ and $d_0$.

\begin{prop}\label{prop-1}
	Assume four given positive real numbers $a_0,b_0,c_0$ and $d_0$ satisfy $a_0+b_0\leq 1$ and $ c_0+d_0\leq 1$.
	
	i) If $a_0/ (a_0+b_0)>c_0/(c_0+d_0)$, let 
	\begin{equation}\label{sim-large}	
\left\{ 
	\begin{split}
	& A_l= \tan( {\arctan({c_0\over d_0}) \over 4}), B_l=\tan({\pi\over 4}+{\arctan({a_0\over b_0})\over 2}),\\
    & a_1={A_l(a_0-B_lb_0)\over A_l-B_l}, b_1={a_0-B_lb_0\over A_l-B_l}, a_2=a_0-a_1, b_2=b_0-b_1;\\
    & C_l= \tan( {\arctan({c_0\over d_0}) \over 2}), D_l=\tan({3\pi\over 8}+{\arctan({a_0\over b_0})\over 4}),\\
    & c_1={C_l(c_0-D_ld_0)\over C_l-D_l}, d_1={c_0-D_ld_0\over C_l-D_l}, c_2=c_0-c_1, d_2=d_0-d_1.
	\end{split}
\right.		
	\end{equation}
	Then, (\ref{decom-1}) and (\ref{decom-2}) hold and $${a_1\over a_1+b_1}<{c_1\over c_1+d_1},\,\ \,\  {a_2\over a_2+b_2}<{c_2\over c_2+d_2}. $$ 
		ii) If $a_0/ (a_0+b_0)<c_0/(c_0+d_0)$,  let 
			\begin{equation}\label{sim-small}	
			\left\{ 
			\begin{split}
				& A_s= \tan( {\arctan({a_0\over b_0}) \over 2}), B_s=\tan({3\pi\over 8}+{\arctan({c_0\over d_0})\over 4}),\\
				& a_1={A_s(a_0-B_sb_0)\over A_s-B_s}, b_1={a_0-B_sb_0\over A_s-B_s}, a_2=a_0-a_1, b_2=b_0-b_1;\\
				& C_s= \tan( {\arctan({a_0\over b_0}) \over 4}), D_s=\tan({\pi\over 4}+{\arctan({c_0\over d_0})\over 2}),\\
				& c_1={C_s(c_0-D_sd_0)\over C_s-D_s}, d_1={c_0-D_sd_0\over C_s-D_s}, c_2=c_0-c_1, d_2=d_0-d_1.
			\end{split}
			\right.		
		\end{equation}
	Then, (\ref{decom-1}) and (\ref{decom-2}) hold and $${a_1\over a_1+b_1}>{c_1\over c_1+d_1}, \,\ \,\ {a_2\over a_2+b_2}>{c_2\over c_2+d_2}. $$ 
\end{prop}

{\bf Proof}. From the top plot in Figure \ref{figure-simpson}-c, $(b_1,a_1)$ is the intersection of the following two lines 
\begin{eqnarray*}
 {y \over x}=A_l,\,\ \,\ {y-a_0\over x-b_0}=B_l.
\end{eqnarray*}
Similarly, from the bottom plot in Figure \ref{figure-simpson}-c, $(d_1,c_1)$ is the intersection of the following two lines 
\begin{eqnarray*}
		 {y \over x}=C_l, \,\ \,\
		 {y-c_0\over x-d_0}=D_l.
\end{eqnarray*}
These justify the choices in (\ref{sim-large}), and the claims in part i) of Proposition \ref{prop-1} (i.e., those below (\ref{sim-large})) follows Lemma~\ref{lem-simpson-1}.

Part ii) of Proposition \ref{prop-1} is established similarly. \fbox{}

Using the Simpson's Paradox as a basic block the general $n$-factor Simpson's Paradox can be built by $2^n-1$ Simpson's Paradoxes. 

\section{An example} For illustration purpose, we now construct a $n$-factor Simpson's Paradox for $n=3$ with the binary variable $X$ and the factor of interest $A$, three additional factors $B_1$, $B_2$ and $B_3$ and the initial values $(a_0,b_0,c_0,d_0)=(0.8,0.2,0.6,0.4)$. 

This 3-factor Simpson's Paradox consists of 7 Simpson's Paradoxes. 
When $B_1$, $B_2$ and $B_3$ are included in the study sequentially, we construct 1 Simpson's Paradox, 2 Simpson's Paradoxes and 4 Simpson's Paradoxes, respectively. Note two facts: i) the Simpson's Paradoxes constructed later do not affect those constructed earlier. i.e., once a Simpson's Paradox is constructed, it is not changed. ii) A new Simpson's Paradox only depends on a single Simpson's Paradox in the previous level. These make the construction similar to building with Lego blocks.

The initial values yield $$	P(X_1|A_1) ={a_0\over a_0+b_0}=0.8> 0.6={c_0\over c_0+d_0}=P(X_1|A_0).$$
For simplicity, we write this relation as $A_1>A_0$ by dropping $X_1$ and $P$.
The first paradox is $$A_1>A_0\,\ \mbox{vs.}\,\ \,\ A_1\bar{B}_1 < A_0\bar{B}_1 \,\ \& \,\ A_1\bar{B}_0 < A_0\bar{B}_0$$
by adding $B_1$.
The second and third paradoxes are
$$A_1\bar{B}_1 < A_0\bar{B}_1 \,\ \mbox{vs.}\,\ A_1\bar{B}_{1,1} > A_0\bar{B}_{1,1} \,\ \& \,\ A_1\bar{B}_{1,0} > A_0\bar{B}_{1,0},$$
$$A_1\bar{B}_0 < A_0\bar{B}_0 \,\ \mbox{vs.}\,\ A_1\bar{B}_{0,1} > A_0\bar{B}_{0,1} \,\ \& \,\ A_1\bar{B}_{0,0} > A_0\bar{B}_{0,0},$$
respectively, by adding $B_2$.
The 4-th through 8-th paradoxes are
$$ A_1\bar{B}_{1,1} > A_0\bar{B}_{1,1} \,\ \mbox{vs.} \,\ A_1\bar{B}_{1,1,1} < A_0\bar{B}_{1,1,1}\,\ \& \,\
A_1\bar{B}_{1,1,0} < A_0\bar{B}_{1,1,0},$$
$$ A_1\bar{B}_{1,0} > A_0\bar{B}_{1,0} \,\ \mbox{vs.} \,\ A_1\bar{B}_{1,0,1} < A_0\bar{B}_{1,0,1}\,\ \& \,\
A_1\bar{B}_{1,0,0} < A_0\bar{B}_{1,0,0},$$
$$ A_1\bar{B}_{0,1} > A_0\bar{B}_{0,1} \,\ \mbox{vs.} \,\ A_1\bar{B}_{0,1,1} < A_0\bar{B}_{0,1,1}\,\ \& \,\
A_1\bar{B}_{0,1,0} < A_0\bar{B}_{0,1,0},$$
$$ A_1\bar{B}_{0,0} > A_0\bar{B}_{0,0} \,\ \mbox{vs.} \,\ A_1\bar{B}_{0,0,1} < A_0\bar{B}_{0,0,1}\,\ \& \,\
A_1\bar{B}_{0,0,0} < A_0\bar{B}_{0,0,0},$$
respectively, by adding $B_3$.

The first paradox is obtained by applying (\ref{sim-large}) to $(a_0^{(1)},b_0^{(1)},c_0^{(1)},d_0^{(1)})=(a_0,b_0,c_0,d_0)$. The superscript ``(i)'' is the index for the $i$-th paradox. Here, $i=1$. We have $(a_1^{(1)},b_1^{(1)},a_2^{(1)},b_2^{(1)})=(0.0263, 0.1047, 0.7737, 0.0953)$ and $(c_1^{(1)},d_1^{(1)},c_2^{(1)},d_2^{(1)})=(0.2010,0.3755, 0.3990, 0.0245)$, indicating
$$P(X_1|A_1B_{1,1})={a_1^{(1)}\over a_1^{(1)}+b_1^{(1)}}=0.2005<0.3486={c_1^{(1)}\over c_1^{(1)}+d_1^{(1)}}=P(X_1|A_0B_{1,1}), $$
$$P(X_1|A_1B_{1,0})={a_2^{(1)}\over a_2^{(1)}+b_2^{(1)}}=0.8904<0.9422={c_2^{(1)}\over c_2^{(1)}+d_2^{(1)}}=P(X_0|A_0B_{1,0}).$$
The second and third paradoxes are constructed  by applying (\ref{sim-small}) to $(a_0^{(2)},b_0^{(2)},c_0^{(2)},d_0^{(2)})=(a_1^{(1)},b_1^{(1)},c_1^{(1)},d_1^{(1)})$ and $(a_0^{(3)},b_0^{(3)},c_0^{(3)},d_0^{(3)})=(a_2^{(1)},b_2^{(1)},c_2^{(1)},d_2^{(1)})$, respectively. Then, we obtain $(a_1^{(2)},b_1^{(2)},a_2^{(2)},b_2^{(2)})$ and $(c_1^{(2)},d_1^{(2)},c_2^{(2)},d_2^{(2)})$ for the second paradox and $(a_1^{(3)},b_1^{(3)},a_2^{(3)},b_2^{(3)})$ and
\begin{table}[H]
	\centering
	\begin{tabular}{l|cccccccc}
		\hline
		&  \multicolumn{8}{c}{\textbf{$A_1$}} \\
		$P(X_1|A_1)$ & \multicolumn{8}{c}{$0.8^{1}$} \\
		&  \multicolumn{4}{c}{\textbf{$A_1\bar{B}_1$}} & \multicolumn{4}{c}{\textbf{$A_1\bar{B}_0$}}\\
		$P(X_1\bar{B}_i|A_1)$ & \multicolumn{4}{c}{0.0263} & \multicolumn{4}{c}{0.7737}\\
		$P(X_1|A_1\bar{B}_i)$ & \multicolumn{4}{c}{$0.2005^{2}$} & \multicolumn{4}{c}{$0.8904^{3}$}\\ 
		&  \multicolumn{2}{c}{\textbf{$A_1\bar{B}_{1,1}$}} 
		& \multicolumn{2}{c}{\textbf{$A_1\bar{B}_{1,0}$}} 
		& \multicolumn{2}{c}{\textbf{$A_1\bar{B}_{0,1}$}} 
		& \multicolumn{2}{c}{\textbf{$A_1\bar{B}_{0,0}$}}\\
		$P(X_1\bar{B}_{i,j}|A_1)$ & \multicolumn{2}{c}{0.0125} & \multicolumn{2}{c}{0.0138}& \multicolumn{2}{c}{0.0748} & \multicolumn{2}{c}{0.6990}\\
		$P(X_1|A_1\bar{B}_{i,j})$ & \multicolumn{2}{c}{$0.1099^{4}$} & \multicolumn{2}{c}{$0.7833^{5}$}& \multicolumn{2}{c}{$0.4693^{6}$} & \multicolumn{2}{c}{$0.9849^{7}$}\\	
		&  \textbf{$A_1\bar{B}_{1,1,1}$} & \textbf{$A_1\bar{B}_{1,1,0}$} 
& \textbf{$A_1\bar{B}_{1,0,1}$} & \textbf{$A_1\bar{B}_{1,0,0}$} 
& \textbf{$A_1\bar{B}_{0,1,1}$} & \textbf{$A_1\bar{B}_{0,1,0}$} 
& \textbf{$A_1\bar{B}_{0,0,1}$} & \textbf{$A_1\bar{B}_{0,0,0}$}\\
$P(X_1\bar{B}_{i,j,k}|A_1)$ & 0.0014 & 0.0111 & 0.0005 & 0.0133& 0.0048 & 0.0700 & 0.0021 & 0.6968\\
$P(X_1|A_1\bar{B}_{i,j,k})$ & 0.0151$^{8}$ & 0.5308$^{9}$& 0.2086$^{10}$ & 0.8805$^{11}$ & 0.0832$^{12}$& 0.6894$^{13}$ & 0.2884$^{14}$ & 0.9924$^{15}$\\	
		\hline
		&  \multicolumn{8}{c}{\textbf{$A_0$}} \\ 
		$P(X_1|A_0)$ & \multicolumn{8}{c}{$0.6^{1}$} \\
		&  \multicolumn{4}{c}{\textbf{$A_0\bar{B}_1$}} & \multicolumn{4}{c}{\textbf{$A_0\bar{B}_0$}}\\
		$P(X_1\bar{B}_0|A_0)$ & \multicolumn{4}{c}{0.2021} & \multicolumn{4}{c}{0.3990}\\
		$P(X_1|A_0\bar{B}_0)$ & \multicolumn{4}{c}{$0.3468^{2}$} & \multicolumn{4}{c}{$0.9422^{3}$}\\
		&  \multicolumn{2}{c}{\textbf{$A_0\bar{B}_{1,1}$}} 
		& \multicolumn{2}{c}{\textbf{$A_0\bar{B}_{1,0}$}} 
		& \multicolumn{2}{c}{\textbf{$A_0\bar{B}_{0,1}$}} 
		& \multicolumn{2}{c}{\textbf{$A_0\bar{B}_{0,0}$}}\\
		$P(X_1\bar{B}_{i,j}|A_0)$ & \multicolumn{2}{c}{0.0163} & \multicolumn{2}{c}{0.1847}& \multicolumn{2}{c}{0.0047} & \multicolumn{2}{c}{0.3943}\\
		$P(X_1|A_0\bar{B}_{i,j})$ & \multicolumn{2}{c}{$0.0579^{4}$} & \multicolumn{2}{c}{$0.6253^{5}$}& \multicolumn{2}{c}{$0.2747^{6}$} & \multicolumn{2}{c}{$0.9703^{7}$}\\	
		&  \textbf{$A_0\bar{B}_{1,1,1}$} & \textbf{$A_0\bar{B}_{1,1,0}$} 
		& \textbf{$A_0\bar{B}_{1,0,1}$} & \textbf{$A_0\bar{B}_{1,0,0}$} 
		& \textbf{$A_0\bar{B}_{0,1,1}$} & \textbf{$A_0\bar{B}_{0,1,0}$} 
		& \textbf{$A_0\bar{B}_{0,0,1}$} & \textbf{$A_0\bar{B}_{0,0,0}$}\\
		$P(X_1\bar{B}_{i,j,k}|A_0)$ & 0.0080 & 0.0082 & 0.0578 & 0.1269 & 0.0022 & 0.0025 & 0.0103 & 0.3840\\
		$P(X_1|A_0\bar{B}_{i,j,k})$ & 0.0298$^{8}$ & 0.7253$^{9}$ & 0.3617$^{10}$ & 0.9367$^{11}$ & 0.1547$^{12}$ & 0.8321$^{13}$ & 0.4923$^{14}$ & 0.9962$^{15}$\\	
	\end{tabular}
\caption{{\it The $3$-factor Simpson's Paradox with factors $A$, $B_1$, $B_2$ and $B_3$. Each pair with the same superscript contains two probabilities for comparison. For example, $0.8^1=P(X_1|A_1)>P(X_1|A_0)=0.6^1$. Each probability is equal to the sum of two probabilities when another $B_i$ is included. For example, $0.0263=P(X_1\bar{B}_1|A_1)= P(X_1\bar{B}_{1,1}|A_1)+P(X_1\bar{B}_{1,0}|A_1)=0.0125+0.0138$.}}
\label{tab-1}
\end{table}
\noindent $(c_1^{(3)},d_1^{(3)},c_2^{(3)},d_2^{(3)})$ for the third paradox.
The 4-th through 7-th paradoxes are constructed similarly. The numerical details of the 7 paradoxes are displayed in Table~\ref{tab-1}. The conclusion about the effect of $A$ on $X$ reverses 3 times as $B_1$, $B_2$ and $B_3$ are included one by one.

\section{Discussion}

Simpson's Paradox is one of the most famous paradoxes in Statistical Science. It essentially says that if we collect data through observational studies, then it is possible that  statistical conclusion on the effect of a certain factor ($A$) on the response ($X$) may reverse as more information in additional factors ($B_i$'s) is collected. 

The classic Simpson's Paradox only involves $X$, $A$ and $B_1$.  The conclusion of the effect of $A$ on $X$ is opposite when $B_1$ is absent or present. We find a connection between this paradox with the plane geometry as in Lemmas ~\ref{lem-simpson-1} and \ref{lem-simpson-1-1}. An R-code is available from the authors to implement the numerical construction in Proposition~\ref{prop-1}. Although we have applied the difference of two proportions, $p_1=P(X_1|A_1\bar{B}_{s_m})$ and $p_0=P(X_1|A_0\bar{B}_{s_m})$, to measure the effect of $A$ on $X$, the conclusion is also true if using the relative risk $p_1/p_0$ and the odds ratio $p_1(1-p_0)/[(1-p_1)p_0]$. This is due to the equivalence of  $p_1-p_0>0$, $p_1/p_0>1$,  and $p_1(1-p_0)/[(1-p_1)p_0]>1$. 

We also investigate whether Simpson's Paradoxes occur multiple times as multiple factors 
$B_i$ are included sequentially. Intuitively, one would expect to make more precise inferences as the information builds up. For Simpson's Paradox, however, we find that it may occur every time we include an additional $B_i$. i.e., the inference about the effect of $A$ on $X$ flips over and over, with no end. This shows that big data may not guarantee correct inferences. 

The results developed so far in this paper are based on the assumption that each of $A$ and $B_i$'s assumes two levels.
This assumption can indeed be loosened to allow for any finite number of levels greater than or equal to two. A case that $A$, $B_1$ and $B_2$ all assume three levels is given in the Supplementary Materials.






\section*{Declaration of conflicting interests}
The author(s) declared no potential conflicts of interest with respect to the research, authorship, and/or publication of this article.



\bibliographystyle{plain}

\begin{thebibliography}{100}
	
	
		\bibitem{Bickel-etal: 1972}	
Bickel, P. J., Hammel, E. A. and  O’Connell, J. W. (1975). Sex bias in graduate admissions: data from Berkeley.
{\it Science} \textbf{187(4175)}: 398–404.

	\bibitem{Blyth: 1972}	
	Blyth, C. R. (1972). On Simpson's Paradox and the sure-thing principle. {\it Journal of the American Statistical Association}. \textbf{67}(338): 364–366.
	

	


	
\bibitem{Charig:1986}
Charig, C. R., Webb, D. R., Payne, S. R., and Wickham, J. E. (1986). Comparison of treatment of renal calculi by open surgery, percutaneous nephrolithotomy, and extracorporeal shockwave lithotripsy. {\it British Medical Journal (Clinical Research Edition)} \textbf{292 (6524)}: 879–882.

	\bibitem{Good:1987}
Good, I. J., Mittal, Y. (1987). The amalgamation and geometry of two-by-two contingency
tables. {\it The Annals of Statistics} 15(2), 694-711.

		\bibitem{Kocik:2001}
	Kocik, J. (2001). Proof without words: Simpson’s paradox. {\it Mathematics Magazine} 74 (5),
	399.
	
		\bibitem{Lindley:1981}	
	Lindley, D. V., Novick, M. R. (1981). The role of exchangeability in inference. {\it The Annals
		of Statistics} 9, 45-58.
	
	\bibitem{Pearson:etal：1899}			
	 Pearson, K., Lee, A., Bramley-Moore, L. (1899). Genetic (reproductive) selection: Inheritance of fertility in man, and of fecundity in thoroughbred racehorses. {\it Philosophical Transactions of the Royal Society A.} \textbf{192}: 257–330. 


	
	\bibitem{Simpson:1951}	
 Simpson, E. H. (1951). The Interpretation of Interaction in Contingency Tables. {\it Journal of the Royal Statistical Society, Series B.} \textbf{13}: 238–241.


	\bibitem{Yule:1903}		 
Yule, G. U. (1903). Notes on the theory of association of attributes in statistics. {\it Biometrika} \textbf{2} (2): 121–134. 
	





	
\end{thebibliography}

\end{document}